\documentclass{amsart}
\usepackage{amsmath,amsfonts,amssymb,stmaryrd}
\usepackage{graphicx,psfrag}

\input xy
\xyoption{all}

\newtheorem{thm}{Theorem}
\newtheorem{rk}{Remark}
\newtheorem{prop}{Proposition}
\newtheorem{clly}{Corollary}
\newtheorem{lemma}{Lemma}
\newtheorem{defi}{Definition}
\newcommand{\pf}{{\flushleft{\bf Proof: }}}

\newcommand{\noi}{\noindent}
\newcommand{\fin}{
			\begin{flushright}
			$\Box$
			\end{flushright}}
\newcommand{\calb}{\mbox{$\mathcal {B}$}}

\newcommand{\ag}{\leftharpoonup}

\begin{document}
\title{Group-Like algebras and Hadamard matrices}

\author{Mariana Haim}
\thanks{The author was partially supported by PEDECIBA and CSIC (Uruguay).}
\date{\today}
\begin{abstract}
We give a description in terms of square matrices of the family of
group-like algebras with $S*id=id*S=u\varepsilon$. In the case that
$S=id$ and $char\Bbbk$ is not $2$ and does not divide the dimension of
the algebra, this translation take us to Hadamard matrices and,
particularly, to examples of biFrobenius algebras satisfying
$S*id=id*S=u\varepsilon$ and that are not Hopf algebras. Finally, we
generalize some known results on separability and coseparability valid
for finite dimensional Hopf algebras to this special class of
biFrobenius algebras with $S*id=id*S=u \varepsilon$, presenting a
version of Maschke's theorem for this family. 
\end{abstract}
\maketitle
\begin{center}
{\large\bf{Introduction.}}\\
\end{center}
\noi
BiFrobenius algebras generalize finite dimensional Hopf algebras, in
the following sense: a biFrobenius algebra has a structure of
Frobenius algebra and coFrobenius coalgebra, linked by a condition
weaker than the pentagonal axiom, namely one assumes the existence of
a (bijective) antimorphism of algebras and coalgebras $S:A \rightarrow
A$. BiFrobenius algebras were introduced in 2000 by Y.Doi and
M. Takeuchi ([DT]).\\ 
\noi
In general, $S$ is not the convolution inverse of the identity. This
is true in the particular situation of Hopf algebras. A natural
question is whether this additional condition implies or not that the
algebra is Hopf; in other words: does the fact that $S$ is the
convolution inverse of the identity forces $A$ to be a bialgebra? \\ 
\noi
In this paper, we show that this is not true by constructing a family
of biFrobenius algebras where $S$ is the convolution inverse of the
identity that are not Hopf algebras.\\ 
Once we know there is an intermediate class between biFrobenius
algebras and finite dimensional Hopf algebras, it could be relevant to
study which properties of finite dimensional Hopf algebras still hold
for this class. In this direction, we give some first steps,
generalizing some known results concerning separability,
semisimplicity and its dual notions.\\ 
\ \\
\noi
In section \ref{uno} we recall the definition and some of the basic
propertis of biFrobenius algebras.\\ 
In section \ref{dos} we present the example of group-like
algebras. Group-like algebras are the natural generalization of group
algebras to the context of biFrobenius algebras. We also study
additional conditions for a group-like algebra to satisfy
$S*id=id*S=u\varepsilon$ and to be a Hopf algebra.\\ 
In section \ref{tres}, we describe, in terms of a family of square
matrices, the algebra structure of a group-like algebra with
$S*id=id*S=u\varepsilon$. As the coalgebra structure of a group-like
algebra is trivial, this matrix approach will give us a complete
description of the family of group-like algebras with
$S*id=id*S=u\varepsilon$.\\ 
In section \ref{cuatro} we briefly present Hadamard matrices. These
matrices will be used in section \ref{cinco} to construct the
counterexamples mentioned above.\\ 
In section \ref{cinco} we treat the cases in which $S=id$ and
$char\Bbbk$ does not divide the dimension of the algebra, showing how
the family of square matrices mentioned before, gives rise, in this
particular case, to a Hadamard matrix. We construct, via Hadamard
matrices, a family, of unbounded dimension, of commutative group-like
algebras with $S*id=id*S=u\varepsilon$ and that are not Hopf
algebras.\\ 
Finally, in section \ref{seis} we prove a version of Maschke's theorem
for this special biFrobenius algebras.\\ 
\begin{section}{BiFrobenius algebras.}\label{uno}
\noi
All along this paper we use Sweedler's notation for the coproduct of a coalgebra.\\
\ \\
Let $A$ be an algebra over $\Bbbk$.\\
We consider $A^*$ with the structure of right $A$-module via $(f \ag
a)(x)=f(ax), \forall f \in A^*, a, x \in A$. For every $f \in A^*$,
the induced map $f \ag :A\rightarrow A^*$ is a morphism of right
$A$-modules.\\ 
If $C$ is a coalgebra over $\Bbbk$, $C$ has a natural structure of
right $C^*$-module, via the right action $c \ag f=\sum f(c_1)c_2,
\forall c \in C, f \in C^*$. For every $c \in C$, the induced map $c
\ag :C^* \rightarrow C$ is a morphism of right $C^*$-modules. 
\begin{defi}\label{f}
Let $\Bbbk$ be a field. A {\bf Frobenius algebra} is a pair $(A,\phi)$
where $A$ is a $\Bbbk$-algebra and  $\phi \in A^*$ is such that the
morphism 
$$
\begin{array}{ll}
\phi \ag :& A \rightarrow A^*,\\
& x \mapsto \phi \ag x
\end{array}
$$
is bijective.\\
\noi
Dually, a {\bf coFrobenius coalgebra} is a pair $(C,t)$ where $C$ is a
$\Bbbk$-coalgebra and $t\in C$ is such that the morphism 
$$
\begin{array}{ll}
t\ag:& C^* \rightarrow C,\\
& f \mapsto t\ag f
\end{array}
$$
is bijective.
\end{defi}
\begin{rk}
\noi
Clearly, Frobenius algebras and coFrobenius coalgebras are finite dimensional.\\
\end{rk}
\noi
We recall now some results valid for the augmented and coaugmented cases.\\
\ \\
\noi
Suppose $A$ is a $\Bbbk$-algebra that admits an algebra morphism
$\varepsilon:A \rightarrow \Bbbk$. In this case, the algebra is said
to be {\bf augmented} and we say that $s \in A$ is a right integral if
$sx=\varepsilon (x)s, \forall x \in A$. 
\begin{prop}\label{int}
If $(A,\phi,\varepsilon)$ is an augmented Frobenius algebra and $s\in
A$ is (the unique element) such that $\phi \ag s=\varepsilon$, then
$s$ is an right integral in $A$. \\ 
\end{prop}
\noi
Dually, suppose $C$ is a $\Bbbk$-coalgebra with a group-like element
$1 \in C$, i.e. $1 \in C$ satisfies $\Delta (1)=1 \otimes 1$. In this
case, we say that the coalgebra is {\bf coaugmented} and that $\psi
\in C^*$ is a right cointegral for $C$ if $\sum \psi(x_1)x_2=\psi(x)1,
\forall x \in C$. 
\begin{prop}\label{coint}
If $(C,t,1)$ is a coaugmented coFrobenius algebra and $\psi \in C^*$
is (the unique element) such that $t\ag \psi=1$, then $\psi$ is a
right cointegral. 
\end{prop}
\noi
The space of right integrals of an augmented Frobenius algebra $A$ is
one dimensional, since  it is $I_r(A)=A^A\cong
(A^*)^A=\Bbbk\varepsilon$.\\ 
\noi
Dually, the space of right cointegrals of a coaugmented coFrobenius
coalgebra is one dimensional.\\ 
\noi
Now we can give the definition of a biFrobenius algebra.
\begin{defi}
Let $\Bbbk$ be a field. The data $(A,m,1,\Delta,\varepsilon,\phi,t)$
is said to be a {\bf biFrobenius} algebra if the following conditions
hold: 
\begin{itemize}
\item $(A,m,1,\phi)$ is a Frobenius algebra,
\item $(A,\Delta,\varepsilon,t)$ is a coFrobenius coalgebra,
\item $\varepsilon$ is an algebra map,
\item $1\in C$ is a group-like element,
\item The map $S:A \rightarrow A$ defined by $S(x)=\sum\phi(t_1x)t_2$
is an antimorphism of algebras and coalgebras. 
\end{itemize}
\end{defi}
\noi
In [DT] it is proved that $S$ is necessarily bijective. $S$ is called
the {\bf antipode} of the biFrobenius algebra $A$.\\ 
\noi
Notice that propositions \ref{int} and \ref{coint} imply that $\phi$ is a cointegral and $t$ is an integral for $A$. Indeed,
$$
S(1)=1, \varepsilon\circ S=\varepsilon,
$$
means that $\sum\phi(t_1)t_2=1$ and that $\phi(tx)=\varepsilon(x), \forall x \in A$;
in other words,
$$
t \ag \phi=1, \mbox{ and } \phi \ag t=\varepsilon.
$$
\begin{defi}
Let $A$ be a biFrobenius algebra, $s,\psi$ respectively a right integral and a right cointegral for $A$. Once can easily show that
$$
s\ag \psi=1 \leftrightarrow \psi(x)=1 \leftrightarrow \psi \ag s=\varepsilon.
$$
We say that $(s,\psi)$ is a biFrobenius pair if one of these (equivalent) conditions holds.
\end{defi}
\noi
Using the fact that the space of right integrals and the space of right cointegrals are one dimensional, it is easy to see that:
\begin{rk}
\ \\
\begin{itemize}
\item Given a non-zero right integral $s\in A$, there is one and only one $\psi \in A^*$ such that $(s,\psi)$ is a biFrobenius pair.
\item Given a non-zero right cointegral $\psi\in A^*$, there is one and only one $s \in A$ such that $(s,\psi)$ is a biFrobenius pair.
\item If $(t,\phi)$ is a biFrobenius pair then any other biFrobenius pair is of the form $(\lambda t, \frac{1}{\lambda}\phi), \lambda \in \Bbbk, \lambda \neq 0$.
\item If $(s,\psi)$ is a biFrobenius pair, then $\sum \psi(s_1x)s_2=S(x)$.
\end{itemize}
\end{rk}
\noi
{\bf\underline{Examples.}}\\
\ \\
\noi
It is well known that every finite dimensional Hopf algebra is biFrobenius. \\
In the next section we present the coalgebraically trivial example of biFrobenius algebras. See [DT], [F] and [WZ] for other examples.\\
\begin{rk}
If $A$ is a biFrobenius algebra, in the finite dimensional algebra $Hom_\Bbbk(A,A)$ the condition $S*id=u\varepsilon$ is equivalent to $id*S=u\varepsilon$.
\end{rk}
\begin{prop}\label{bis}
Let $A$ be a biFrobenius algebra over a field $\Bbbk$. If $A$ is also a bialgebra then $S*id=u\varepsilon$. Hence $A$ is a Hopf algebra.
\end{prop}
\pf
As $A$ is a bialgebra if and only if $\Delta$ is an algebra morphism, then
$$
\begin{array}{ll}
(S*id)(x)&=\sum S(x_1)x_2=\sum \phi(t_1x_1)t_2x_2=\sum \phi\left((tx)_1\right )(tx)_2=\\
&=\phi(tx)1=(\phi \ag t)(x)1=\varepsilon(x)1.
\end{array}
$$
\fin
\noi
We will show in section \ref{cuatro} that the converse is not true, namely that we can find biFrobenius algebras with the property $S*id=u\varepsilon$ that are not Hopf algebras, in the sense that the multiplication and the comultiplication are not compatible.
\begin{rk}\label{debil}
Notice that, in order to prove proposition \ref{bis}, we have used in fact a condition that is weaker than the compatibility between product and coproduct in a bialgebra. Namely the condition
$$
\varepsilon(x)\sum t_1 \otimes t_2=\sum(tx)_1 \otimes (tx)_2=\sum t_1x_1 \otimes t_2x_2 \ \ \ \ \ (*)
$$
implies that $S*id=u\varepsilon$. It is relevant to know that condition $(*)$ is in fact weaker than the compatibility condition for a bialgebra. See the end of section \ref{cinco}.
\end{rk}
\end{section}
\begin{section}{Group-Like Algebras.}\label{dos}
\noi
We present here a family of examples of biFrobenius algebras, introduced in [D] by Y.Doi. It generalizes the group algebras in the sense that the coalgebraic structure of the group-like algebras is trivial, as it is for group algebras.\\
\noi
We will present this example by adding, step by step, the minimal conditions required to have a biFrobenius structure, coalgebraically trivial, in a finite dimensional $k$-vector space.\\
\ \\
\noi
All along the rest of the paper we set
$$
I=\{0,1,2,\cdots,n\}.
$$
\noi
We will use $1$ for both $1_\Bbbk$ and $1_A$, since it will be clear from the context to which of them we are referring.\\
\ \\
\noi
We start by giving  the minimal data we need in order to have the simplest coFrobenius coalgebra structure in a finite dimensional vector space. We leave the proof of the next proposition to the reader.
\begin{prop}\label{ev}
Let $\Bbbk$ be a field and $A$ be a finite dimensional $\Bbbk$-vector space with basis $\calb =\{b_0,b_1,b_2,\cdots b_n\}$. Let $\varepsilon:\calb \rightarrow \Bbbk$ be such that $\varepsilon(b_i)\neq 0, \forall i\in I$.
\begin{enumerate}
\item If we define
$$
\Delta(b_i)=\frac{1}{\varepsilon(b_i)}b_i \otimes b_i, \forall i\in I, \ \ \ \ \ t=b_0+b_1+b_2\cdots+b_n
$$
and we extend $\Delta$ and $\varepsilon$ to $A$ by linearity, then $(A,\Delta,\varepsilon,t)$ is a coFrobenius coalgebra.
\item Assume moreover that $\varepsilon(b_0)=1$. Then, with respect to
the structure $\Delta$ given above, $b_0$ is a group-like element of
$A$ and the preimage $\phi \in A^*$ of $b_0$ by $t\ag$ satisfies
$\phi(b_i)=\delta_{i,0}, \forall i \in I$. 
\end{enumerate}
\end{prop}
\noi
Assume moreover that we have an algebra structure in the vector space $A$ given, in terms of the structure constants, by
$$
b_ib_j=\sum_{k\in I}p_{ij}^kb_k, \forall i,j \in I, \ \ \ \ 1=b_0.
$$
\noi
Notice that the fact that $b_0=1$ forces the cointegral $\phi \in A^*$ to be as defined in proposition \ref{ev}, (2).
\noi
Next proposition presents necessary and sufficient conditions on the coefficients $p_{ij}^k$ in order to make the map $S:A\rightarrow A$ defined by the formula $S(x)=\sum\phi(t_1x)t_2$ an antimorphism of algebras and of coalgebras. Notice that
$$
S(b_j)=\sum_{i\in I}\frac{1}{\varepsilon(b_i)}\phi(b_ib_j)b_i=\sum_{i\in I}\frac{1}{\varepsilon(b_i)}p_{ij}^0b_i.
$$
\noi
If $\sigma:I \rightarrow I$ is a permutation, we say that $S$ extends $\sigma$ if $S(b_i)=b_{\sigma(i)}, \forall i \in I$.
\begin{prop}\label{anti}
Let $\Bbbk$ be a field and $A$ be a finite dimensional $\Bbbk$-algebra
with basis $\calb =\{b_0,b_1,\cdots,b_n\}$ and $b_0=1$. Let $p_{ij}^k,
i,j,k \in I$ be the structure constants associated to the algebra $A$
in the basis $\calb$. Take $\varepsilon:\calb \rightarrow \Bbbk$ to be
a function that is never $0$ and $\sigma:I\rightarrow I$ a permutation
of $I$ and endow $A$ with the comultiplication defined before. Then:\\

\begin{enumerate}
\item \label{a} $S$ extends $\sigma$ if and only if $p_{ij}^0=\varepsilon(b_i)\delta_{i,\sigma(j)}$.
\item Suppose \em({\ref{a})} holds. Then:
\begin{enumerate}
\item $S$ is an antimorphism of coalgebras if and only if $\varepsilon (b_{\sigma(i)})=\varepsilon(b_i), \forall i \in I.$
\item $S$ is an antimorphism of algebras if and only if
$$\sigma(0)=0, \mbox{ and }p_{ij}^k=p_{\sigma(j)\sigma(i)}^{\sigma(k)}.$$
\end{enumerate}
\item If $S$ is an antimorphism of algebras and coalgebras that extends $\sigma$, then $\sigma^2=id$.
\end{enumerate}
\end{prop}
\pf
\begin{enumerate}
\item We know that $S(b_j)=\sum_{i\in
I}\frac{1}{\varepsilon(b_i)}p_{ij}^0 b_i$. As $\calb$ is linearly
independent, the condition $S(b_j)=b_{\sigma(j)}, \forall j \in I$ is
satisfied if and only if
$\frac{1}{\varepsilon(b_i)}p_{ij}^0=\delta_{i,\sigma(j)}$, which is
equivalent to $p_{ij}^0=\varepsilon(b_i)\delta_{i,\sigma(j)}$. 
\item We assume now that $S(b_i)=b_{\sigma(i)}, \forall i \in I$.
\begin{enumerate}
\item We have that $\frac{1}{\varepsilon (b_{\sigma(i)})}
b_{\sigma(i)}\otimes
b_{\sigma(i)}=\Delta(b_{\sigma(i)})=\Delta(S(b_i))$ and $(S \otimes
S)(\Delta(b_i))=\frac{1}{\varepsilon(b_i)}b_{\sigma(i)}\otimes
b_{\sigma(i)}$. Then $S$ is anticommutes with $\Delta$ if and only if
$\varepsilon (b_{\sigma(i)})=\varepsilon(b_i), \forall i \in
I$. Notice that this implies that $S$ is counital , i.e. that
$\varepsilon \circ S=\varepsilon$. 
\item $S$ is an antimorphism of algebras if and only if $S(1)=1$ and $S(xy)=S(y)S(x), \forall x,y \in A$.
But $S(1)=1$ if and only if $S(b_0)=b_0$, which means, in terms of $\sigma$, that $b_{\sigma(0)}=b_0$, i.e. $\sigma(0)=0.$\\
On the other hand $S(xy)=S(y)S(x), \forall x,y \in A$ if and only if
$S(b_ib_j)=S(b_j)S(b_i), \forall i,j \in I$. Writing this equality in
terms of the structure constants, we obtain $S\left(\sum_{k\in
I}p_{ij}^k b_k\right)=b_{\sigma(j)}b_{\sigma(i)}, \forall i,j \in I$
and then $ \sum_{k\in I} p_{ij}^k b_{\sigma(k)}=\sum_{k\in I}
p_{\sigma(j)\sigma(i)}^kb_k, \forall i,j \in I$. Then $\sum_{k\in I}
p_{ij}^k b_{\sigma(k)}=\sum_{k\in I}
p_{\sigma(j)\sigma(i)}^{\sigma(k)}b_{\sigma(k)}, \forall i,j \in I$
i.e. $p_{ij}^k=p_{\sigma(j)\sigma(i)}^{\sigma(k)}, \forall i,j,k \in
I$. \\ 
\end{enumerate}
\item Assume that $S$ is an antimorphism of algebras and coalgebras that extends $\sigma$.
Observe that, using (1) and (2), we have that $\forall i \in I:
\varepsilon\left
(b_{\sigma(i)}\right)=\varepsilon(b_{\sigma(i)})\delta_{\sigma(i),\sigma(i)}=p_{\sigma(i)i}^0=p_{\sigma(i)\sigma^2(i)}^{\sigma(0)}=p_{\sigma(i)\sigma^2(i)}^0=\varepsilon\left
(b_{\sigma(i)}\right ) \delta_{\sigma(i)\sigma^3(i)}$, hence
$\sigma^3(i)=\sigma(i),\forall i \in I$ and, as $\sigma $ is
bijective, $\sigma^2(i)=i, \forall i \in I$. 
\end{enumerate}
\fin
\noi
Now we are able to summarize all the conditions required in order to
have in the vector space $A$ a structure of biFrobenius algebra. The
following theorem will follow directly from propositions \ref{ev} and
\ref{anti}. 
\begin{thm}\label{gl}
Let $\Bbbk$ be a field and $A$ be a finite dimensional $\Bbbk$-vector
space with basis $\calb =\{b_0, b_1, \cdots, b_n\}$. Let $p: A\otimes
A\rightarrow A$ be a linear map with structure constants
$\{p_{ij}^k\mid i,j,k \in I\}$ with respect to $\calb$. Let
$\varepsilon:\calb\rightarrow \Bbbk$ be a function and $\sigma:I
\rightarrow I$ be a permutation, where $I=\{0,1,2, \cdots n\}$. If we
have that 
$$
\begin{array}{ll}
(GL1) \ p \mbox{ is associative}, &(GL2) \ p \mbox{ has neutral element }b_0=1,\\
(GL3) \ \varepsilon(b_i)\neq 0, \forall i \in I,& (GL4) \ \varepsilon(b_{\sigma(i)})=\varepsilon(b_i), \forall i \in I,\\
(GL5) \ p_{ij}^0=\delta_{i,\sigma(j)}\varepsilon(b_i),&(GL6) \ p_{ij}^k=p_{\sigma(j)\sigma(i)}^{\sigma(k)},\\
(GL7) \ \sigma(0)=0, &(GL8) \ \varepsilon \mbox{ extends to a morphism of algebras }\varepsilon:A \rightarrow \Bbbk,
\end{array}
$$
then $A$ has a structure of biFrobenius algebra with $\Delta, \phi$ and $t$ defined as before.
\end{thm}
\pf
Notice that the conditions $(GL1)$ to $(GL8)$ guarantee, following
propositions \ref{ev} and \ref{anti} that $A$ is a Frobenius algebra,
$A$ is a coFrobenius coalgebra and that the induced map $S$ is an
antimorphism of algebras and coalgebras. Notice that $\varepsilon:A
\rightarrow \Bbbk$ is a morphism of algebras (condition
$(GL8)$). Moreover, if $(GL8)$ holds, then $\varepsilon(1)=1$, hence
$\Delta (1)=\Delta(b_0)=\frac{1}{\varepsilon(b_0)}b_0\otimes
b_0=1\otimes 1$ and then $1$ is a group-like element for $A$. 
\fin
Now we present the definition of a group-like algebra (see [D]).
\begin{defi}
A $\Bbbk$-{\bf group-like algebra} is a 5-tuple $(A,\calb
,p,\varepsilon,\sigma)$ such that $A$ is a $\Bbbk$-algebra with basis
$\calb =\{b_0,b_1,\cdots b_n\}$, $p:A\otimes A \rightarrow A$ is a
linear map, $\varepsilon:\calb \rightarrow k$ is a function and
$\sigma:I \rightarrow I$ is a bijection, such that conditions $(GL1)$
to $(GL8)$ are satisfied. 
\end{defi}
\begin{clly}
Every group-like algebra is a cocommutative biFrobenius algebra with $S^2$=id.
\end{clly}
\noi
Next we will exhibit the additional conditions needed in order to have a Hopf algebra structure in $A$.
\begin{prop}\label{grupo}
Let ${\bf A}=(A,\calb ,p,\varepsilon,\sigma)$ be a $\Bbbk$-group-like algebra. Then:
\begin{enumerate}
\item $S * id=u\varepsilon$ if and only if $\calb $ is closed under the operation of taking inverses with respect to $p$.
\item The following conditions for ${\bf A}$ are equivalent:
\begin{enumerate}
\item [(i)]  ${\bf A}$ is a Hopf algebra,
\item [(ii)] ${\bf A}$ is a bialgebra,
\item [(iii)]$(\calb ,p,1)$ is a monoid,
\item [(iv)] $(\calb ,p,1)$ is a group.
\end{enumerate}
\item In any of the previous situations we have that $\varepsilon(b_i)=1, \forall i\in I.$
\end{enumerate}
\end{prop}
\pf
\begin{enumerate}
\item
If $S*id=u\varepsilon$, then $(S*id)(b_i)=\varepsilon(b_i)b_0, \forall i \in I$. Then we have
$$
\frac{1}{\varepsilon(b_i)}b_{\sigma(i)}b_i=\varepsilon(b_i)b_0,
$$
and
$\frac{1}{\varepsilon(b_i)}p_{\sigma(i)i}^0=\varepsilon(b_i)$. Then
$p_{\sigma(i)i}^0=\left(\varepsilon(b_i)\right)^2$. But conditions
(GL5) and (GL4) give
$p_{\sigma(i)i}^0=\varepsilon(b_{\sigma(i)})=\epsilon(b_i)$, so using
condition (GL3) we get that $\forall i \in I: \varepsilon(b_i)=1$. 
Then we have that $b_{\sigma(i)}b_i=b_0=1$ (and also that $b_ib_{\sigma(i)}=b_0$, since $A$ is finite dimensional), which means $b_i^{-1}=b_{\sigma(i)}$.\\
For the converse direction, assume $b_i$ has an inverse in
$\calb$. From the condition (GL5) we deduce that the inverse of $b_i$
has to be $b_{\sigma(i)}$. From the equality $b_{\sigma(i)}b_i=b_0$
and conditions (GL8) and (GL4) we get
$\varepsilon(b_i)^2=\varepsilon(b_0)=1$, hence
$\varepsilon(b_i)\in\{1.-1\}$, which implies
$\frac{1}{\varepsilon(b_i)}=\varepsilon(b_i)$. We get that
$\frac{1}{\varepsilon(b_i)}b_{\sigma(i)}b_i=\varepsilon(b_i)b_0=\varepsilon(b_i)1,
\forall i \in I$ and hence, the equality $S*id=u\varepsilon$ holds in
the basis $\calb$ and then in ${\bf A}$. 
\item The fact that ${\bf A}$ is a Hopf algebra if and only if it is a
bialgebra, follows directly from proposition \ref{bis}. Hence $(i)
\leftrightarrow (ii)$. 
Suppose now that ${\bf A}$ is a bialgebra. We know, by proposition
\ref{bis} that $S*id=u\varepsilon$. So, by the proof of (1), we get
that $\varepsilon(b_i)=1, \forall i \in I$.\\ 
As ${\bf A}$ is a bialgebra, $\forall i,j \in I$, we have that
$
\Delta(b_ib_j)=\Delta(b_i)\Delta(b_j),
$
which means, in terms of the structure constants,
$$
\sum_{k\in I}p_{ij}^k (b_k \otimes b_k)=\sum_{k,l \in I}p_{ij}^kp_{ij}^l(b_k \otimes b_l),\forall i,j \in I.
$$
Then, $p_{ij}^kp_{ij}^l=0$ for each pair $(k,l)$ such that $k \neq l$,
which means that given $i,j\in I$ there is $k_0 \in I$ such that
$p_{ij}^{k} = 0, \forall k\neq k_0.$ Moreover, using the bialgebra
condition above; for $k=l=k_0$, we get that $p_{ij}^{k_0}=\left
(p_{ij}^{k_0}\right)^2$, so that $p_{ij}^{k_0}\in \{0,1\}$. If
$p_{ij}^{k_0}=0$, we would have $b_i.b_j=0$ and so, using condition
(GL8), $\varepsilon(b_i)\varepsilon(b_j)=0$ which contradicts
condition (GL3).\\ 
We conclude that $p_{ij}^{k_0}=1$ and then $b_i.b_j=b_{k_0}$.\\
Conversely, suppose $\calb$ is closed under $p$. Then for each pair
$(i,j) \in I \times I$ there is $k\in I$ such that
$b_i.b_j=b_k$. Applying $\Delta$ we obtain that 
$$
\begin{array}{ll}
\Delta(b_ib_j)&=\Delta(b_k)=\frac{1}{\varepsilon(b_k)}b_k\otimes b_k=\frac{1}{\varepsilon(b_i)\varepsilon(b_j)}(b_ib_j\otimes b_ib_j)=\\
&=\left ( \frac{1}{\varepsilon(b_i)} b_i \otimes b_i\right)\left (\frac{1}{\varepsilon(b_j)} b_j \otimes b_j\right)=\Delta(b_i)\Delta(b_j)
\end{array}
$$
and then we have proved $(ii)\leftrightarrow (iii)$. (Notice we have used condition (GL8)).
Finally, as ${\bf A}$ is a Hopf algebra if and only if ${\bf A}$ is a
bialgebra and $S*id=u\varepsilon$, and these two conditions hold if
and only if $(\calb ,p,1)$ is a monoid closed under taking inverses
i.e. if $(\calb, p, 1)$ is a group, we have proved $(i)\leftrightarrow
(iv)$. 
\item
It is clear from the proof of (1) and (2) that in both cases $\varepsilon(b_i)=1, \forall i \in I$.
\end{enumerate}
\fin
\noi
In particular, we have the following result.
\begin{clly}\label{cor1}
If a group-like algebra ${\bf A}=(A,\calb,p,\varepsilon,\sigma)$ is
Hopf, then $(\calb,p,1)$ is a group, $\sigma$ is the inverse in
$\calb$ and ${\bf A}$ is the group algebra $\Bbbk\calb$. 
\end{clly}
\end{section}
\begin{section}{The matricial approach.}\label{tres}
\noi
In the next theorem we give a description in terms of square matrices
of the family of group-like algebras with the property that $S*id =
u\varepsilon$. In order to do this, we need some previous
definitions. From now on, we denote by $w$ both a vector in
$\Bbbk^{n+1}$ and its transpose. Also $\{e_0,e_1,\cdots,e_n\}$ is the
canonical basis of $\Bbbk^{n+1}$. 
\begin{defi}\label{sist}
Let $\sigma: I \rightarrow I$ be a permutation and $v =(1, 1, . .,.1) \in \Bbbk^{n+1}$.
\begin{itemize}
\item A $(I, \sigma)$-{\bf system} is a family of invertible matrices
$\mathcal{F} = \{F_i\mid i\in I\}\subseteq M_{n+1}(\Bbbk)$ satisfying:

$$
F_{\sigma(i)}= F^t_i= F^{-1}_i, \forall i \in I,\ \ \    F_iv = v, \forall i \in I, \ \ \    F_ie_0 = e_i, \forall i \in I.
$$
\item If $\mathcal{F}$ and $\mathcal{G}$ are two systems, we say that they are {\bf compatible} if $F_iG_j = G_jF_i, \forall i, j \in I.$
\end{itemize}
\end{defi}\label{lio}
\begin{lemma}\label{lema}
If $\mathcal{F}$ and $\mathcal {G}$ are compatible $(I, \sigma)$-systems, then
$$
F_0=G_0=Id \mbox{  and  }
(F_i)_{kj}= (G_j)_{ki}.
$$
\end{lemma}
\pf
Observe that $F_0e_i=F_0G_ie_0=G_iF_0e_0=G_ie_0=e_i$.\\
On the other hand, we have that $F_ie_j = F_iG_je_0 = G_jF_ie_0 = G_je_i$. Therefore
$$
(F_i)_{kj}= (F_ie_j)_k = (G_je_i)_k = (G_j)_{ki}
$$
and the proof is finished.
\fin
\noi
Notice that we can consider, without loss of generality, group-like
algebras whose supporting vector space is $\Bbbk ^{n+1}$. We will do
so from now on, letting $I=\{0,1,2,\cdots,n\}$,
$\mathcal{C}=\{e_0,e_1,\cdots e_n\}$ be the canoncial basis of
$\Bbbk^{n+1}$ and for each $i, j,k \in I, p_{ij}^k\in \Bbbk$. \\ 
What follows is the framework needed in order to formulate next theorem. We define
\begin{itemize}
\item $p:\Bbbk ^{n+1}\otimes \Bbbk ^{n+1} \rightarrow \Bbbk ^{n+1}, \ p(e_i \otimes e_j)=\sum_{k\in I} p_{ij}^ke_k$,
\item ${\bf 1}: \Bbbk ^{n+1} \rightarrow \Bbbk, \ {\bf 1}(e_i)=1,
\forall i \in I$, 
\item $\mathcal{H}=\{H_i\mid i \in I\}$, where $H_i\in M_{n+1}(\Bbbk)$ is defined by $(H_i)_{kj}= p_{ij}^k, \forall i,j,k \in I$,
\item $\mathcal{V}=\{V_i\mid i \in I\}$, where $V_i\in M_{n+1}(\Bbbk)$ is defined by $(V_j)_{ki}=p_{ij}^k, \forall i, j,k \in I$.
\end{itemize}
\begin{rk}\label{asoc}
Notice that:
$$
\left\{
\begin{array}{l}
H_i \mbox{ is the matrix associated to the map } e_i.\_: A \rightarrow
A \mbox{ in the canonical basis},\\ 
V_j \mbox{ is the matrix associated to the map } \_.e_j : A
\rightarrow A \mbox{ in the canonical basis}. 
\end{array}
\right.
$$
\end{rk}
\begin{thm}\label{teo2}
In the context defined above, if $\sigma: I \rightarrow I$ is an
arbitrary permutation, we have that $(k^{n+1},\mathcal{C},p,{\bf
1},\sigma)$ is a group-like algebra with $S * id =u\varepsilon$ if and
only if $\mathcal{H}$ and $\mathcal{V}$ are compatible
$(I,\sigma)$-systems and $\sigma^2=id$. 
\end{thm}
\pf
Assume first that $(\Bbbk ^{n+1},\mathcal {C} , p,{\bf 1},\sigma)$ is
a group-like algebra. We know that $\sigma^2=id$ (corollary
\ref{cor1}). By remark \ref{asoc}, it is clear that $H_ie_0 = V_ie_0 =
e_i, \forall i \in I$ and also that $H_{\sigma(i)}= H^{-1}_i$ and
$V_{\sigma(i)}= V^{-1}_i$ (since $p$ is associative and
$S*id=u\varepsilon$ implies that $e_{\sigma(i)}e_i= e_0 =
e_ie_{\sigma(i)}, \forall i \in I$).\\ 
As we have that
$$
\phi\left (e_{\sigma(k)}(e_ie_j) \right)=\phi\left ((e_{\sigma(k)}e_i)e_j\right),, \forall, i,j,k \in I,
$$
we deduce that $\sum_{l\in I}p_{\sigma(k)l}^0p_{ij}^l=\sum_{l\in
I}p_{\sigma(k)i}^lp_{lj}^0$, which implies, by condition $(GL5)$ (and
the fact that $\sigma^2=id$), that 
$$
p_{ij}^{k}=p_{\sigma(k)i}^{\sigma(j)}, \forall i,j,k \in I.
$$
Using condition $(GL6)$, we get that $p_{ij}^k=p_{\sigma(i)k}^j,
\forall i,j,k \in I$, which means that
$H_i^t=H_{\sigma(i)}$. Similarly, we prove that
$V_i^t=V_{\sigma(i)}$.\\ 
\ \\
Now,
$$
H_iv =\left (\sum_{j \in I}p_{ij}^0,\sum_{j \in
I}p_{ij}^1,\cdots,\sum_{j\in I}p_{ij}^n\right)=\left (\sum_{j \in
I}p_{\sigma(i)0}^j,\sum_{j \in I}p_{\sigma(i)1}^j,\cdots,\sum_{j\in
I}p_{\sigma(i)n}^j\right)=v, 
$$
where last equality follows from
$\varepsilon(b_{\sigma(i)})b_0)=\varepsilon(b_{\sigma(i)}\varepsilon(b_0)=1$
(condition $(GL8)$). In a similar way we prove that $V_iv = v, \forall
i \in I$.\\ 
We have already proved that the families $\mathcal{H}$ and $\mathcal{V}$ are $(I, \sigma)$-systems.
We still have to prove that they are compatible.
For this, we use the associativity of $p$, i.e.  $(e_ie_k)e_j = e_i(e_ke_j)$.  If we put this equality in terms of the associated matrices we get
$$
H_iV_j=V_jH_i.
$$
and we are done.\\
\ \\
Let us prove the converse, so assume that $\mathcal{H}$ and
$\mathcal{V}$ are $(I, \sigma)$-compatible systems.  We have to check
conditions (GL1) to (GL8) of theorem \ref{gl} and also that $S * id =
u\varepsilon$.\\ 
The associativity follows from $V_jH_i= H_iV_j, \forall i, j \in I$ so we have (GL1).\\
By lemma \ref{lema}, we have that $H_0 = V_0 = Id$ and this implies,
by remark \ref{asoc}, that $e_0 = 1$ is the neutral element of the
algebra, i.e. (GL2).\\ 
Conditions (GL3) and (GL4) hold trivially, since $\varepsilon ={\bf 1}$.\\
As $H_{\sigma(i)}=H_i^t$, we have that
$$
p^0_{ij}= (H_i)_{0j}=\left (H_{\sigma(i)}\right)_{j0}=p_{\sigma(i)0}^j= (V_0)_{j\sigma(i)}= \delta_{\sigma(i),j},
$$
so condition (GL5) holds (notice we have used $\sigma^2=id$ and $\varepsilon = {\bf 1}$).\\
Now, using first that $H_j^t= H_{\sigma(j)}$, then that $V_k^t=V_{\sigma(k)}$, and finally that $H_i^t=H_{\sigma(i)}$, we deduce that
$$
p^{\sigma(k)}_{\sigma(j)\sigma(i)}= p^{\sigma(i)}_{j\sigma(k)} = p^j_{\sigma(i)k}= p_{ij}^k,
$$
i.e. condition (GL6).\\
Let us prove condition (GL7), i.e. $\sigma(0) = 0$. We have that
$H_{\sigma(0)}=H_0^t=Id$, but notice that $H_i \neq Id, \forall i \neq
0$, since $H_ie_0 = e_i$, so we get $\sigma(0) = 0$.\\ 
To prove that $\varepsilon = {\bf 1}$ is a morphism of algebras, it is
enough to verify that $\forall i, j \in I$,
$\varepsilon(e_ie_j)=\varepsilon(e_i)\varepsilon(e_j)$. But
$H_{\sigma(i)}= H^t_i$and $H_iv = v, \forall i \in I$ imply that 
$$
\sum_{k \in I}p_{ij}^k= \sum_{k \in I}p_{\sigma(i)k}^j=1, \forall i, j \in I,
$$
i.e. $\varepsilon(e_ie_j) = \sum_{k \in I}p_{ij}^k= 1 = \varepsilon(e_i)\varepsilon(e_j)$.\\
It remains to prove that $S * id = u\varepsilon$.  This can be deduced directly from associativity and the fact that $H_{\sigma(i)}= H^{-1}_i$.
\fin
\begin{clly}\label{cor2}
Let $(I, \sigma)$ be as in definition \ref{sist} and $\mathcal{C} $ be the canonical basis of $\Bbbk ^{n+1}$.
\begin{enumerate}
\item There is a one to one correspondence between group-like algebras
of the form ${\bf A} = (\Bbbk ^{n+1}, \mathcal{C}, p, \varepsilon,
\sigma)$ satisfying $S * id = u\varepsilon$ and pairs of compatible
$(I,\sigma)$-systems. 
\item Let ${\bf A}$ be as described in (1) and
$(\mathcal{H},\mathcal{V})$ its corresponding pair of compatible
$(I,\sigma)$-systems. Then {\bf A} is a Hopf algebra if and only if
the family $\mathcal{H}$ is a group with the usual matrix product (and
this happens if and only if the family $\mathcal{V}$ is a group with
the usual matrix product). 
\end{enumerate}
\end{clly}
\pf
\begin{enumerate}
\item We recall that in proposition \ref{grupo} we proved that in
group-like algebras with $S * id = u\varepsilon$ it holds that
$\varepsilon ={\bf 1}$. It follows from theorem \ref{teo2} that we can
define correspondences 
$$
\begin{array}{lll}
(\Bbbk ^{n+1}, \mathcal{C}, p, {\bf 1},\sigma) &\mapsto &(\mathcal{H}, \mathcal{V} )\mbox{ with }(H_i)_{kj}= p_{ij}^k= (V_j)_{ki},\\
\ \\
(\mathcal{H}, \mathcal{V} )  &\mapsto&  (\Bbbk ^{n+1}, \mathcal{C}, p,{\bf 1}, \sigma), \mbox{ where}\\
&&                      p : \Bbbk ^{n+1}  \otimes \Bbbk ^{n+1} \rightarrow \Bbbk ^{n+1} \mbox{ is given by } p_{ij}^k= (H_i)_{kj}.\\
&& (\mbox{its structure constants with respect to }\mathcal{C}).
\end{array}
$$
\noi
Lemma \ref{lema} guarantees that these correspondences induce a bijection (more specifically, that the first map is surjective).
\item By proposition \ref{grupo}, ${\bf A}$ is a Hopf algebra if and
only if $\calb$ is closed under $p$ if and only if $\forall i,j \in
I$, there is $k \in I$ such that $e_ie_j = e_k$ if and only if
$\forall i, j \in I$, there is $k \in I$ such that $H_iH_j = H_k$ (or
$V_iV_j = V_k)$ and we are done. 
\end{enumerate}
\fin
\end{section}
\begin{section}{Hadamard matrices}\label{cuatro}
\noi
The Hadamard maximal determinant problem is the following: find the
matrices of a given size with entries $+1$ and $-1$ with the largest
(in absolute value) possible determinant. Despite well over a century
of work, beginning with Sylvester's results of 1867, the general
question remains unanswered. However, it is known that Hadamard
matrices are solutions to Hadamard's maximal $n^2$-determinant
problem. 
\begin{defi}\label{had}
A Hadamard matrix of size $n$ is a square matrix (of size $n \in
\mathbb N$) with coefficients in $\{1, -1\}$ that is orthogonal (with
the usual inner product of $\mathbb C^n$) and such that the first row
and the first column are both the vector $v = (1, 1,\cdots,1)$. 
\end{defi}
\begin{rk}
Sometimes the condition concerning the first column and the first row
(called the normalization condition) is not required in the definition
of a Hadamard matrix. However, if we multiply by $-1$ any row or any
column of an orthogonal matrix with coefficients in $\{-1,1\}$, we
still get an orthogonal matrix. In this sense, any "generalized"
Hadamard matrix is equivalent to a normalized one. 
\end{rk}
\noi
It is easy to check that if ${\bf P}$ is a Hadamard matrix of size $n$
then $n = 1, 2$ or a positive multiple of $4$. Also, it is clear that
defining 
$$
{\bf P}_1 = (1),\ \
{\bf P}_2 =
\left(\begin{array}{rr}
1&1\\ 1&-1
\end{array}\right),\ \
{\bf P}_k =\left(\begin{array}{rr}   {\bf P}_{k-1}&{\bf P}_{k-1}\\{\bf P}_{k-1}&-{\bf P}_{k-1}\end{array}\right)
$$
we obtain Hadamard matrices of size $2^k$.\\
For the purposes of the applications to group-like algebras, it will
be relevant whether there exist or not Hadamard matrices whose size is
not a power of $2$. The answer is negative: it has been proved that
there are Hadamard matrices of unbounded size $n$ where $n$ is not a
power of $2$. 
Indeed, Paley's construction (see [W])  guarantees that there exist
Hadamard matrices of size $n$, for every $n$ divisible by $4$ and of
the form $2^e(p^m + 1)$, with $m, e, p \in\mathbb N, m \neq 0, p
\mbox{ prime }, p \neq 2$. We give below an example of a Hadamard
matrix of size $12$ (corresponding to $ p = 5, e = m = 1$).\\ 
$$
\left(
\begin{array}{rrrrrrrrrrrr}
1 & 1&  1 & 1& 1 & 1&1 & 1& 1& 1& 1& 1\\
\ \\
1 &{\bf -1}& {\bf -1} & 1&{\bf -1} &{\bf -1}&{\bf -1}& 1& 1& 1&{\bf -1}& 1\\
\ \\
1 & 1& {\bf -1} &{\bf -1}& 1 &{\bf -1}&{\bf -1}&{\bf -1}& 1& 1& 1&{\bf -1}\\
\ \\
1 &{\bf -1}&  1 &{\bf -1}&{\bf -1} & 1&{\bf -1}&{\bf -1}&{\bf -1}& 1& 1& 1\\
\ \\
1 & 1& {\bf -1} & 1&{\bf -1} &{\bf -1}& 1&{\bf -1}&{\bf -1}&{\bf -1}& 1& 1\\
\ \\
1 & 1&  1 &{\bf -1}& 1 &{\bf -1}&{\bf -1}& 1&{\bf -1}&{\bf -1}&{\bf -1}& 1\\
\ \\
1 & 1&  1 & 1&{\bf -1} & 1&{\bf -1}&{\bf -1}& 1&{\bf -1}&{\bf -1}&{\bf -1}\\
\ \\
1 &{\bf -1}&  1 & 1& 1 &{\bf -1}& 1&{\bf -1}&{\bf -1}& 1&{\bf -1}&{\bf -1}\\
\ \\
1 &{\bf -1}& {\bf -1} & 1& 1 & 1&{\bf -1}& 1&{\bf -1}&{\bf -1}& 1&{\bf -1}\\
\ \\
1 &{\bf -1}& {\bf -1} &{\bf -1}& 1 & 1& 1&{\bf -1}& 1&{\bf -1}&{\bf -1}& 1\\
\ \\
1 & 1& {\bf -1} &{\bf -1}&{\bf -1} & 1& 1& 1&{\bf -1}& 1&{\bf -1}&{\bf -1}\\
\ \\
1 &{\bf -1}&  1 &{\bf -1}&{\bf -1} &{\bf -1}& 1& 1& 1&{\bf -1}& 1&{\bf -1}
\end{array}
\right)
$$
\end{section}
\begin{section}{The case $S = id$.}\label{cinco}
\noi
In this section, we consider the particular case of group-like
algebras supported in $\Bbbk ^{n+1}$ with $S * id =u\varepsilon$,
where $S = id$. They are of the form 
$$
{\bf A}=(\Bbbk ^{n+1}, \mathcal{C} , p, {\bf 1}, id),\mbox{ satisfying also that } e^2_i= 1, \forall i \in I.
$$
\noi
This special type of group-like algebras are commutative, since $id : {\bf A} \rightarrow {\bf A}$ is an antimorphism of algebras.
\ \\
\noi
This section is divided in three main parts:
\begin{itemize}
\item First we consider the correspondence described in corollary
\ref{cor2} for this particular case.  In other words, we describe all
possible pairs of $(I, id)$-compatible systems. 
\item Second, we show that this family of all pairs of $(I,
id)$-compatible systems is in one to one correspondence with the
family of Hadamard matrices. 
\item Finally, we use some results on Hadamard matrices (mentioned in
section 4) to find examples of group-like algebras satisfying $S * id
= u\varepsilon$ that are not Hopf algebras. 
\end{itemize}
\noi
Suppose we have a group-like algebra ${\bf A} = (\Bbbk ^{n+1},
\mathcal{C}, p, {\bf 1}, id)$ with $e^2_i= 1, \forall i \in I$ and
consider the bijection described in corollary \ref{cor2}. 
Observe that, as ${\bf A}$ is commutative, we have that $p_{ij}^k=
p_{ji}^k, \forall i,j,k \in I$ and therefore 
$$
H_i= V_i, \forall i \in I,
$$
meaning that the compatible pair associated to ${\bf A}$ is $(\mathcal{H},\mathcal{H})$, so it can be thought
as a single $(I, id)$-system, with the additional condition that
$$
H_iH_j =H_jH_i, \forall i, j \in I.
$$
This can be summarized as follows:
\begin{defi}
Let $\mathcal{H}$ be a $(I, id)$-system. We say that $\mathcal{H}$ is
{\bf self-compatible} if $H_iH_j = H_jH_i, \forall i, j \in
I$. Explicitely the set $\{H_0,H_1,\cdots,H_n\}$ satisfies 
$$
H_i=H_i^{-1}=H_i^t,\ \ \ \ H_iH_j=H_jH_i,\ \ \ \ H_iv=v, \ \ \ \ H_ie_0=e_i,
$$
where $i,j \in I, v=(1,1,\cdots,1)$.
\end{defi}
\begin{thm}
Let $\Bbbk$ be a field. There is a one to one correspodence between
group-like algebras of the form $(\Bbbk ^{n+1}, \mathcal{C}, p,{\bf
1}, id)$ satisfying also that $e_i^2= 1, \forall i \in I $ and
self-compatible $(I, id)$-systems. 
\end{thm}
\noi
Next we describe, provided some minor restrictions on $char\Bbbk$, a
bijective correspondence between $(I, id)$-self-compatible systems and
Hadamard matrices.\\ 
In the rest of the paper all our systems will be taken over $(I, id)$.
\begin{thm}\label{teo4}
Let $\Bbbk$ be a field and $n\in \mathbb N$. Then:
\begin{enumerate}
\item Assume that $char \Bbbk \neq 2$. Then the matrices of a
self-compatible system $\mathcal{H}=\{H_0,H_1,\cdots,H_n\}$ admit a
basis of common eigenvectors $\{v^0, v^1, . .v^n\}$ such that the
matrix 
$$
{\bf P} =\left (v^0|v^1|\cdots|v^n\right)
$$
is a Hadamard matrix.\\
Moreover, we have that
$$
H_i={\bf P}Di{\bf P}^{-1},
$$
where $D_i$ is a diagonal matrix whose diagonal is the $i^{\mbox{th}}$-row of ${\bf P}$. (Explicitly $(D_i)_{rs}= \delta_{r,s}{\bf P}_{ir}$.)
\item Conversely, assume now that $char \Bbbk$ does not divide
$n+1$. Then, given a Hadamard matrix {\bf P}, we can construct a
self-compatible system $\mathcal{H} = \{H_i \mid i \in I\}$ such that
the columns of ${\bf P}$ form a common basis of eigenvectors for the
matrices of $\mathcal{H}$. 
\end{enumerate}
\end{thm}
\pf
\begin{enumerate}
\item First notice that, as $H_i^2=Id$, each $H_i$ is equivalent to a
diagonal matrix with diagonal entries in $\{1,-1\}$. Denote as
$\langle, \rangle$ the standardt inner product in $\Bbbk ^{n+1}$.We
have also that $H_iH_j = H_jH_i, \forall i, j \in I$ and that
$H_i^t=H_i$, i.e. each $H_i, i \in I$ is symmetric. Therefore, we can
obtain for the matrices in $\mathcal{H}$ a common basis
$\calb=\{v^0,v^1,\cdots,v^n\}$ of orthogonal eigenvectors.\\ 
As $v = (1, 1, \cdots,1)$ is fixed for all $H_i, i \in I$, we can choose $v^0=v$.\\
\ \\
Let $w = (w_0, w_1, \cdots w_n) \in \Bbbk^{n+1}$ be any common non
zero eigenvector.\\ 
Then, for each $i \in I, H_iw = w$ or $H_iw = -w.$ But $H_iw = \sum_{k \in I} w_k(H_ie_k)$
and $(H_ie_k)_i = \langle e_i,H_ie_k\rangle = \langle
H_i^te_i,e_k\rangle = \langle H^{-1}_ie_i,e_k\rangle =\langle e_0,
e_k\rangle = \delta_{0,k}$ 
hence $(H_iw)_i= w_0$ and therefore, as $(H_iw)_i=\pm w_i$, 
$$
 \forall i \in I : w_i= w_0 \mbox{ or }w_i= -w_0.
$$
Take $w$ such that $w_0=1$. Writing $H_iw = \lambda_iw$ and computing
the $i$-th coordinate we deduce that $1=w_0=(H_iw)_i=\lambda_iw_i$ and
therefore $\lambda_i=w_i$. In other words, for each $i$ of $H_i$, the
eigenvalue of $H_i$ associated to $w$ is $w_i$, i.e. $H_iw=w_iw$.\\ 
We can assume $(v^j)_0=1$ and still have an orthogonal basis of eigenvectors. Call ${\bf P}$ the matrix 
$$
{\bf P}=\left(v^0 |v^1 | \cdots|v^n\right).
$$
Then ${\bf P}$ has coefficients in $\{1, -1\}$ and
$$
H_i={\bf P}D_i{\bf P}^{-1},
$$ 
where $D_i$ is a diagonal matrix formed by the (ordered) eigenvalues
of $H_i$. As the eigenvalue associated to $v^j$ for $H_i$ is
$(v^j)_i$, we have that the $i^{\mbox{th}}$-row of ${\bf P}$ gives all
the eigenvalues, meaning the principal diagonal in $D_i$ is the
$i^{\mbox{th}}$-row of ${\bf P}$.\\ 
Moreover, ${\bf P}$ is an orthogonal matrix with coefficients in
$\{1,-1\}$  whose first row and first column are $v = (1, 1, . .,.1)$,
i.e. ${\bf P}$ is a Hadamard matrix.\\  
\item Take a Hadamard matrix ${\bf P}$ of size $n + 1$.\\
It is clear that ${\bf P}$ is invertible since ${\bf P}{\bf
P}^t=(n+1)Id$ and $char \Bbbk$ does not divide $n+1$. Let $r_i$ be the
$i^{\mbox{th}}$-row of ${\bf P}$.  Consider $D_i$ the diagonal matrix
where the principal diagonal is $r_i$ and take 
$$
H_i= {\bf P}D_i{\bf P}^{-1}.
$$
We have to show that $\mathcal{H} = \{H_i\mid i \in I\}$ is a self-compatible system. As $ H_i =\frac{1}{n+1}{\bf P}D_i{\bf P}^t$, then
$$
\begin{array}{ll}
(H_i)^t  &= H_i,\\
H_iH^t_i &= H^2_i= \frac{1}{(n+1)^2}{\bf P}D_i{\bf P}^t{\bf P}D_i{\bf P}^t=\frac{1}{n+1}{\bf P}(D_i)^2{\bf P}^t=\frac{1}{n+1}{\bf P}{\bf P}^t=Id
\end{array}
$$
Now, as ${\bf P}e_0={\bf P}^t e_0=v$ and ${\bf P}^tv=(n+1)e_0$, we have that, $\forall i \in I$, 
$$
\begin{array}{ll}
H_iv = \frac{1}{n+1}{\bf P}D_i{\bf P}^tv=\frac{1}{n+1}{\bf P}D_i(n+1)e_0={\bf P}e_0=v,\\
H_ie_0=\frac{1}{n+1}{\bf P}D_i{\bf P}^te_0=\frac{1}{n+1}{\bf P}D_iv=\frac{1}{n+1}{\bf P}r_i=e_i.
\end{array}
$$
We have proved that $\mathcal{H}$ is a system. It remains to prove
that it is self-compatible, i.e. that the matrices in $\mathcal{H}$
commute with each other, but this follows from the fact that they are
equivalent to a diagonal matrix via the same matrix ${\bf P}$. 
\end{enumerate}
\fin
\noi
For the next theorem we will consider in $\Bbbk ^{n+1}$ the so called point-wise product
$$
w.w'= (w_0w'_0, w_1w'_1, \cdots,w_nw'_n).
$$
\noi
\begin{rk}\label{prod}
If ${\bf P}$ is a Hadamard matrix, the rows in ${\bf P}$ are self-inverses under the point-wise product.\\
On the other hand, the matrices in a $(I,id)$-compatible system are self-inverses under the usual matrix product.\\
Next theorem gives a stronger link between these two products.
\end{rk}
\begin{thm}\label{teo5}
Let $\mathcal{H}$ be a self-compatible system and ${\bf P}$ a Hadamard
matrix of common eigenvectors for $\mathcal{H}$. Then the matrices in
$\mathcal{H}$ form a group (with the usual matrix product) if and only

if the rows in ${\bf P}$ form a group (with the point-wise product).
\end{thm}
\pf
\noi
In view of remark \ref{prod}, we only have to prove that $\mathcal{H}$ is closed under the usual
matrix product if and only if the set of rows in ${\bf P}$ is closed under the point-wise
product. But this follows directly from
$$
H_iH_j = H_k \mbox{ if and only if }D_iD_j = D_k \mbox{ if and only if } r_ir_j = r_k,
$$
where $r_i$ is the $i^{\mbox{th}}$-row of ${\bf P}$.
\fin
\begin{rk}
The existence of an invertible Hadamard matrix of size $n+1$ implies
that $char \Bbbk$ does not divide $n+1$ (since the matrix is
invertible) and in particular that $char \Bbbk \neq 2$ (since the
existence of a Hadamard matrix implies that $n+1$ is even). 
\end{rk}
\begin{clly}\label{ce}
\noi
There are biFrobenius algebras of unbounded dimension satisfying $S * id = u\varepsilon$ that are not Hopf algebras.
\end{clly}
\pf
\noi
Take $\Bbbk$ such that $char \Bbbk =0$. Let $k \in \mathbb N$. As we
observed before, it is known that there is a Hadamard matrix of size
$n + 1$ bigger 
than $k$, with $n + 1$ not a power of $2$.\\
Let $\mathcal{H}$ be the self-compatible system constructed from $
{\bf P}$ as in theorem \ref{teo4} and let 
{\bf A} be its associated group-like algebra as in corollary
\ref{cor2}. We claim that ${\bf A}$ is not a Hopf algebra.\\ 
If ${\bf A}$ were Hopf algebra, then, by corollary \ref{cor2},
$\mathcal{H}$ is a group with the usual matrix- 
product and therefore, by theorem \ref{teo5}, the rows of ${\bf P}$ would form a group. Then we would have an
abelian finite group all whose elements haver order two and whose size
is not a power of $2$, and this contradicts the structure theorem of
finite abelian groups. Therefore ${\bf A}$ is not a Hopf algebra. 
\fin
\noi
We finish this section by proving that the condition $(*)$ we
considered in remark \ref{debil} is weaker than the compatibility
condition for bialgebras. 
\begin{prop}\label{t}
Let ${\bf A}=(A,\calb,p,{\bf 1},\sigma)$ be a group-like algebra and $t\neq 0$ be a right integral for ${\bf A}$. Then $S*id=u\varepsilon$ if and only if
$$
\forall x \in A, \epsilon(x)t=\Delta(tx)=\sum t_1x_1\otimes t_2x_2. \ \ \ \ \ \ (*)
$$
\end{prop}
\pf
The converse implication was proved in proposition \ref{bis}.\\
For the direct implication, notice that, as $S*id=u\varepsilon$, we
have $\varepsilon={\bf 1}$, and then $tb_j=t, \forall j \in
I$. Therefore $\Delta(tb_j)=\Delta(t), \forall j \in I$ and condition
$(*)$ becomes $\Delta(t)=\sum t_1b_j\otimes t_2b_j, \forall j \in I$,
i.e. 
$\sum_{k \in I} b_k \otimes b_k=\sum_{i,k,l \in I} p_{ij}^kp_{ij}^l b_k\otimes b_l, \forall j \in I.$
In other words, we have to prove that
$$
\sum_{i \in I}p_{ij}^kp_{ij}^l=\delta_{k,l}, \forall j,k,l \in I.
$$
But $S*id=u\varepsilon$ implies $b_{\sigma(j)}(b_jb_{\sigma(l)})=b_{\sigma(l)}, \forall j,l \in I$ and therefore
$$
\sum_{i,u\in I}p_{\sigma(j)i}^up_{j\sigma(l)}^ib_u=b_{\sigma(l)}, \forall i, j,l \in I.
$$
For $u=\sigma(k)$ we have $\sum_{i\in I}p_{\sigma(j)i}^{\sigma(k)}p_{j\sigma(l)}^i=\delta_{\sigma(k),\sigma(l)}$ or equivalently,
$$
\sum_{i \in I}p_{\sigma(j)\sigma(i)}^{\sigma(k)}p_{j\sigma(l)}^{\sigma(i)}=\delta_{k,l}.
$$
But $p_{\sigma(j)\sigma(i)}^{\sigma(k)}=p_{ij}^{k}$ (from condition
(GL6)) and $p_{j\sigma(l)}^{\sigma(i)}=p_{l\sigma(j)}^i=p_{ij}^l$
(from condition (GL6) and from the fact that $V_{\sigma(j)}=V_j^t$).\\

Therefore, we have that $\sum_{i \in I}p_{ij}^kp_{ij}^l=\delta_{k,l}$ and we are done.
\fin
The following result follows directly from corollary \ref{ce} and proposition \ref{t}.
\begin{clly}
There are biFrobenius algebras of unbounded dimension satisfying $\Delta(tx)=\sum t_1x_1 \otimes t_2x_2$ for a right integral $t \neq 0$ and that are not bialgebras.
\end{clly}
\end{section}
\begin{section}{On separability and coseparability of biFrobenius algebras.}\label{seis}
\noi
The notion of separability is classical in ring theory. Every
separable algebra is semisimple. Indeed, a $\Bbbk$-algebra $A$ is
separable if and only if for any field extension $E \supseteq \Bbbk$,
the $E$-algebra $A^E=A\otimes_{\Bbbk} E$ is semisimple. \\ 
We use an alternative definition -see Definition \ref{sep}- of a
separable algebra and give also the (dual) notion of a coseparable
coalgebra (see [DMI] and [T]).\\ 
\ \\
\noi
In this section, if $A$ is an algebra and $C$ is a coalgebra, we consider
\begin{itemize}
\item $A$ with the usual structure of $A-A$ bimodule,
\item $A\otimes A$ with the structure of $A-A$ bimodule given by $a (x\otimes y)b=ax\otimes yb$,
\item $C$ with the usual structure of $C-C$ bicomodule,
\item $C \otimes C$ with the structure of $C-C$ bicomodule given by $\sum (x\otimes y)_{-1} \otimes (x\otimes y)_0 \otimes (x\otimes y)_1=\sum x_1 \otimes x_2 \otimes y_1 \otimes y_2$.
\end{itemize}
It is evident that the product $m:A\otimes A \rightarrow A$ is an
epimorphism of $A-A$ bimodules and that the coproduct $\Delta: C
\rightarrow C \otimes C$ is a monomorphism of $C-C$ bicomodules. 
\begin{defi}\label{sep}
\begin{enumerate}
\item A $\Bbbk$-algebra $A$ is said to be {\bf separable} if the product $m:A \otimes A \rightarrow A$ splits in the category of $A-A$ bimodules.
\item A $\Bbbk$-coalgebra $C$ is said to be {\bf coseparable} if the
coproduct $\Delta:C \rightarrow C \otimes C$ splits in the category of
$C-C$ bicomodules. 
\end{enumerate}
\end{defi}
\noi
The following result is well known and can be found for example in [S].
\begin{thm}
Let $H$ be a finite dimensional Hopf algebra and $t\in H$ be a
non-zero right or left integral. Then the following statements are
equivalent. 
\begin{enumerate}
\item[(i)] $H$ is semisimple.
\item[(ii)] $H$ is separable.
\item[(iii)] $\varepsilon (t)\neq 0$.
\end{enumerate}
\end{thm}
\noi
We give a proof of this result in the context of biFrobenius algebras
satisfying $S *id=u\varepsilon$. We refer to [D] for a slightly
different proof of $(c) \rightarrow (b)$. We also state its dual
version. 
\begin{thm}
Let $(A,m,1,\Delta,\varepsilon,t,\phi,S)$ be a biFrobenius algebra such that $S*id=id*S=u\varepsilon$. Then:
\begin{enumerate}
\item The following assertions are equivalent:
\begin{enumerate}
\item[(a)] $A$ is semisimple.
\item[(b)] $A$ is separable.
\item[(c)] $\varepsilon(t)\neq 0$.
\end{enumerate}
\item The following assertions are equivalent:
\begin{enumerate}
\item[(a')] $A$ is cosemisimple.
\item[(b')] $A$ is coseparable.
\item[(c')] $\phi(1)\neq 0$.
\end{enumerate}
\end{enumerate}
\end{thm}
\pf
\begin{enumerate}
\item
\noi
For (a) implies (c) take $ker(\varepsilon)\subseteq A$ as a right
$A$-submodule. As $A$ is semisimple, there is a right submodule $I
\subseteq A$ of dimension $1$ such that $ker(\varepsilon) \oplus
I=A$. We write $1=a+s$, with $\varepsilon(a)=0$ and $s \in I$, then
$I=\Bbbk s$ and $\varepsilon(s)=1$. If $x \in A$, we have $sx \in I$
and then  $sx=\lambda s$, for some $\lambda \in \Bbbk$. Applying
$\varepsilon$ we deduce that $\varepsilon(x)=\lambda$ and therefore
$sx=\varepsilon(x)s$, which means that $s$ is a right integral for
$A$. Hence, $t=\mu s \in I$, for some $\mu \in \Bbbk, \mu \neq 0$,
which implies $\varepsilon(t)\neq 0$.\\ 
For (c) implies (b), take $\bar S$ the composite inverse of the
identity, $s=\frac{1}{\varepsilon(t)}t$ and $\psi \in A^*$ such that
$(s,\psi)$ is a biFrobenius pair. Consider $\delta:A \rightarrow A
\otimes A, \delta(x)=\sum x\bar S(s_2)\otimes s_1$. From the equality
$\sum x \bar S(s_2)s_1=x$, it follows that $\delta$ splits $m$. We
have to prove that $\delta$ is a morphism of $A-A$ bimodules; in other
words that $\sum axb\bar S(s_2)\otimes s_1=\sum ax\bar S (s_2)\otimes
s_1b$. In order to do this, it is enough to show that 
$$
\sum b\bar S(s_2)\otimes s_1=\sum \bar S(s_2) \otimes s_1b,
$$
which is equivalent, by applying the morphism $id \otimes \psi \ag: A
\otimes A \rightarrow A \otimes A^*$  to the equality $b \bar
S(s_2)\psi(s_1z)=\bar S(s_2) \psi(s_1bz), \forall z \in A,$ or $b\bar
S(Sz)=\bar S(S(bz)), \forall z \in A$. This last equality is obviously
true.\\ 
Finally, let us prove that (b) implies (a). Take $\delta$ the map that
splits the product and put $\delta(1)=\sum_{i=1}^l r_i \otimes
r^i$. We have $\sum r_ir^i=1$ and $\sum ar_i \otimes
r^i=a\delta(1)=\delta(a)=\delta(1)a=\sum r_i \otimes r^ia$.\\ 
Given $N \subseteq M$ an inclusion of $A$-modules, consider a
$\Bbbk$-linear map $p:M \rightarrow N$ such that $p_{|N}=id_N$. The
map 
$$
\begin{array}{ll}
\pi:M \rightarrow N,\\
\pi(x)=\sum_{i=1}^l p(xr_i)r^i.
\end{array}
$$
is a morphism of $A$-modules that splits the projection. Indeed, for $x \in N$, $xr_i \in N, \forall i=1,2,\cdots,l$ and therefore $\pi(x)=\sum xr_ir^i=x\sum r_ir^i=x$.\\
Moreover $\pi$ is a morphism of right $A$-modules, since
$$
\pi(xa)=\sum p(xar_i)r^i=\sum p(xr_i)r^ia=\pi(x)a.
$$
\item The proof of this equivalence is obtained by dualizing the methods described above. We only present a sketch in order to obtain explicit expressions for the morphisms involved.\\
The proof that (a') implies (c') is obtained from $(1)$, working in
$A^*$.\\ 
In the proof that (c') implies (b'), in order to split $\Delta$ in the
category of $A-A$ bicomodules, we take the map $*: A \otimes A
\rightarrow A$ given by $x*y=\sum\psi(Sxy_1)y_2$, where $\psi \in A^*$
is a right cointegral such that $\psi(1)=1$.\\ 
Finally, to prove that (b') implies (c') we take a map $*: A \otimes A
\rightarrow A$ that splits $\Delta$. If $N \subseteq M$ is a right
$A$-subcomodule and $p:M \rightarrow N$ is a linear map that splits
the inclusion, we consider $\pi:M \rightarrow N$, $\pi(x)=\sum
p(m_0)_0\otimes \varepsilon\left (p(m_0)_1\otimes m_1\right)$. 
\end{enumerate}
\fin
For the particular case of group-like algebras, the following version of Maschke's theorem can be easily deduced.
\begin{clly}
Every group-like algebra is coseparable.\\
A group-like algebra is separable if and only if $\sum_{i=0}^n \varepsilon(b_i) \neq 0$.\\
In particular a group-like algebra with $S*id=u\varepsilon$ is separable if and only if $char\Bbbk$ does not divide its dimension.
\end{clly}
\end{section}

\ \\
\ \\
\ \\
\ \\

\footnotesize
\noindent Mariana Haim,\\
{\sc Facultad de Ciencias, Universidad de la Rep\'ublica.}\\ 
{\sc Igu\'a 4225, 11400 Montevideo, Uruguay.}
\begin{flushleft}
{\em E-mail adress}: {\tt negra@cmat.edu.uy}
\end{flushleft}
\end{document}